\documentclass[letterpaper, 10 pt, conference]{ieeeconf}
\IEEEoverridecommandlockouts 
\let\labelindent\relax

\usepackage{amsmath}
\usepackage{amssymb}
\usepackage{enumitem}
\usepackage[colorlinks = true, linkcolor = blue, citecolor = blue]{hyperref}
\usepackage[nameinlink]{cleveref}
\usepackage{graphicx}

\newtheorem{assumption}{Assumption}
\newtheorem{theorem}{Theorem}
\newtheorem{lemma}{Lemma}
\newtheorem{algorithm}{Algorithm}
\newtheorem{definition}{Definition}
\newtheorem{example}{Example}
\newtheorem{remark}{Remark}

\crefname{equation}{}{}
\crefname{section}{Section}{Sections}

\begin{document}

\title{Flexible-step MPC for Unknown Linear Time-Invariant Systems}

\author{Markus Pietschner, \quad Christian Ebenbauer, \quad Bahman Gharesifard, \quad Raik Suttner%
\thanks{Markus Pietschner, Christian Ebenbauer, and Raik Suttner are with the Chair of Intelligent Control Systems at RWTH Aachen University, Aachen, Germany.
        {\tt\small \{markus.pietschner,christian.ebenbauer, raik.suttner\}@ic.rwth-aachen.de}}%
\thanks{Bahman Gharesifard is with the Department of Mathematics and Statistics at Queen's University, Kingston, Canada.
        {\tt\small bahman.gharesifard @queensu.ca}}%
}

\maketitle
\begin{abstract}
We propose a novel flexible-step model predictive control algorithm for unknown linear time-invariant discrete-time systems. The goal is to asymptotically stabilize the system without relying on a pre-collected dataset that describes its behavior in advance. In particular, we aim to avoid a potentially harmful initial open-loop exploration phase for identification, since full identification is often not necessary for stabilization. Instead, the proposed control scheme explores and learns the unknown system online through measurements of inputs and states. The measurement results are used to update the prediction model in the finite-horizon optimal control problem. If the current prediction model yields an infeasible optimal control problem, then persistently exciting inputs are applied until feasibility is reestablished. The proposed flexible-step approach allows for a flexible number of implemented optimal input values in each iteration, which is beneficial for simultaneous exploration and exploitation. A generalized control Lyapunov function is included into the constraints of the optimal control problem to enforce stability. This way, the problem of optimization is decoupled from the problem of stabilization. For an asymptotically stabilizable unknown control system, we prove that the proposed flexible-step algorithm can lead to global convergence of the system state to the origin.
\end{abstract}

\section{Introduction}\label{sec:1}
Optimization-based control methods using receding horizon information, such as model predictive control (MPC), form an important class of control strategies, which are intensively studied from both theoretical and practical point of views \cite{GruneBook}. These methods offer many desirable features; for example, the ability to optimize performance and the capability to include constraints into the underlying control problems. However, a major limitation, which can reduce the practical applicability, is that these methods rely on an accurate prediction model for the dynamics of the to-be-controlled system \cite{Hewing20202}. The derivation of such a prediction model is typically based on physical laws and identification. This can lead to uncertainties in the prediction model, which in turn may have a negative impact on the performance of the control system. In the present paper, we focus on the problem of stabilizing a system by MPC without prior model knowledge.

The problem of dealing with uncertainties in the prediction model is intensively studied in the literature on MPC for unknown systems. A natural approach is to employ tools from the field of \emph{adaptive control}. For example, simultaneous identification and MPC for a class of nonlinear systems is presented in \cite{Marafioti2014}. Dual adaptive MPC for unknown linear systems is studied in \cite{Heirung2017}, where a Kalman filter is used for parameter estimation. However, the above studies do not provide guarantees of closed-loop stability. By contrast, the adaptive MPC scheme in \cite{Goncalves2016} ensures closed-loop stability and constraint satisfaction, which, however, involves a computationally expensive min-max MPC formulation. Under the assumption of knowledge of a nominal linear system as well as bounds on modeling errors, the results in \cite{Aswani2013} establish state and input constraint satisfaction as well as robustness with respect to modeling errors. Based on set-membership identification, the adaptive approach in \cite{Tanaskovic2014} is computational tractable, allows for an online model adaption and possesses provable closed-loop properties, but it is limited to open-loop stable systems. The adaptive MPC scheme in \cite{Kim2008} guarantees feasibility and stability of the closed-loop system, but the approach relies on the conservative assumption that a pre-stabilizing feedback controller is known.

There exists a variety of \emph{data-driven} MPC schemes for unknown systems. Many of these schemes rely on pre-exploration of the unknown system through \emph{persistently exciting} (PE) inputs. It is shown in \cite{WaardeBook} that if pre-exploration provides a sufficiently informative dataset of inputs and states, then it is even possible to compute a stabilizing state-feedback controller for the unknown system. However, such a linear state-feedback approach does not allow for performance optimization and an inclusion of state and input constraints. The studies \cite{Coulson2019} and \cite{Coulson20192} introduce \emph{data-enabled predictive control} (DeePC) for unknown deterministic and stochastic linear time-invariant (LTI) systems. An extension of DeePC to linear time-varying systems with permanent updates of the PE data set is reported in \cite{Baros2022}. The data-driven MPC scheme in \cite{Berberich20202} can even lead to exponential stability under the assumption of a controllable and observable LTI system. The above data-driven MPC methods require access to a pre-recorded data set, which can be obtained by measuring the system's response to a sufficiently long sequence of PE inputs. However, such pre-exploration can cause an undesirable initial behavior of the system, which we want to avoid in the present paper.

Finally, we emphasize that there exist many other strategies for stabilizing unknown systems via MPC. These include, among others, schemes based on reinforcement learning \cite{Berkenkamp2017}, neural networks \cite{Patan2015}, and Gaussian process prediction models \cite{Maiworm2018}. An algorithm based on convex body chasing is proposed in \cite{Yu20232}, which guarantees bounded states in the presence of bounded adversarial disturbances, under the assumption that the unknown system belongs to a \emph{known} convex and compact set of stabilizable systems.

In this paper, we propose a novel approach to unknown systems, which belongs to the class of so-called \emph{flexible-step MPC}. A characteristic feature of flexible-step MPC is that a flexible number of elements of an optimal input sequence is applied to the control system; as for instance in \cite{Yang1993}, \cite{Alamir2017}, \cite{Worthmann2017}, \cite{Furnsinn2024}. Here, we focus on the recently introduced flexible-step approach from \cite{Furnsinn2024}, which is particularly suited for the purpose of stabilization. In particular, this approach has been successfully applied to stabilize nonholonomic systems \cite{Furnsinn2023} as well as switched linear systems \cite{Furnsinn2025}. So far the flexible-step framework from \cite{Furnsinn2024} is only investigated for systems for which an accurate prediction model is known. In the present study, we present a first extension to a class of unknown control systems.

\emph{Contributions:}
We present a new flexible-step MPC approach for unknown LTI systems, which does not require a potentially harmful pre-exploration phase. In particular, we do not aim at full identification of the system. The proposed MPC approach combines exploration and exploitation in the spirit of dual control \cite{Feldbaum1960}. A key feature is the concept of a so-called \emph{generalized discrete-time control Lyapunov function} from \cite{Furnsinn2024}, which is incorporated into the constraints of the optimal optimal control problem. By doing so, the tasks of optimization and stabilization are effectively decoupled. Input and state constraints can be taken into account as soft constraints through the cost function. The proposed MPC method neither requires a zero-state terminal constraint nor a suitable design of the terminal cost (which is generally difficult for an unknown system). Under the natural assumption that the unknown LTI system is stabilizable, our main theoretical result provides sufficient conditions for global convergence of the system state to the origin. Our numerical tests further demonstrate that the flexible-step approach has some inherent robustness properties with respect to measurement noise.

The paper is organized as follows. In \Cref{sec:2}, we provide a precise problem statement, introduce frequently used notation, and explain the main components of flexible-step MPC. The rather short \Cref{sec:3} is devoted to flexible-step MPC in case of a known LTI system for which we can prove global exponential stability. Then, in \Cref{sec:4}, we turn to the case of an unknown LTI system. This section contains the novel flexible-step MPC scheme and the main convergence result. The proposed method is tested through numerical simulations in \Cref{sec:5}. Some conclusions are drawn in \Cref{sec:6}.

\section{Problem statement and main components of flexible-step MPC}\label{sec:2}
We consider a discrete-time linear time-invariant control system with state space $\mathbb{X}$ and input space $\mathbb{U}$, both of which are finite-dimensional vector spaces. The time domain is the set $\mathbb{N}$ of nonnegative integers. Let $\mathcal{L}(\mathbb{X},\mathbb{X})$ (resp.~$\mathcal{L}(\mathbb{U},\mathbb{X})$) denote the space of linear maps from $\mathbb{X}$ to $\mathbb{X}$ (resp.~from $\mathbb{U}$ to $\mathbb{X}$). The state dynamics are determined by two linear and time-invariant maps $A\in\mathcal{L}(\mathbb{X},\mathbb{X})$ and $B\in\mathcal{L}(\mathbb{U},\mathbb{X})$ as follows. If, at time $t\in\mathbb{N}$, the system is in state $x(t)\in\mathbb{X}$ and an input $u(t)\in\mathbb{U}$ is applied, then the state $x(t+1)$ at time $t+1$ is given by
\begin{equation}\label{eq:unknownSystem}
x(t+1) \ = \ A x(t) + B u(t).
\end{equation}
Our goal is to asymptotically stabilize \cref{eq:unknownSystem} about the origin without requiring explicit knowledge of $A$ and $B$. To this end, we will propose a novel flexible-step MPC approach, which can lead to the desired convergence of the system state under the natural assumption that \cref{eq:unknownSystem} is asymptotically stabilizable.

\emph{Notation:}
To distinguish actual states and inputs of \cref{eq:unknownSystem} from predicted states and inputs, we will use the following notational conventions. For state and input trajectories of \cref{eq:unknownSystem}, we use the bracket notation $x(t)$ and $u(t)$ for the values of $x$ and $u$ at $t\in\mathbb{N}$, where $t$ is a typical symbol for the current time of \cref{eq:unknownSystem}. Typical symbols for the time variable of predicted states and inputs are $i$, $k$, and $l$, which are attached as subscripts, e.g., $x_k$ and $u_k$. To describe time domains of predicted states and inputs, for all $i,l\in\mathbb{N}$, we let $[\vphantom{]}i\colon\!l\vphantom{(})$ (resp.~$[i\colon\!l]$) denote the set of $k\in\mathbb{N}$ for which $i\leq{k}<l$ (resp.~$i\leq{k}\leq{l}$). For all sets $I,S$, we let $S^I$ denote set of all maps from $I$ to $S$. In particular, for every $k\in\mathbb{N}$ and every set $S$, we may identify $S^{[\vphantom{]}0\colon\!k\vphantom{(})}$ with the $k$-fold Cartesian product $S^k$, where the values of each $s\in{S^{[\vphantom{]}0\colon\!k\vphantom{(})}}$ are indexed by the set $[\vphantom{]}0\colon\!k\vphantom{(})$; that is, $s=(s_0,\ldots,s_{k-1})$. For all $i,k,l\in\mathbb{N}$ with $i\leq{k}\leq{l}$, every set $S$, and every $s\in{S}^{[\vphantom{]}0\colon\!l\vphantom{(})}$, we let $s_{[\vphantom{]}i\colon\!k\vphantom{(})}$ and $s_{[i\colon\!k]}$ denote the subsequences $(s_i,\ldots,s_{k-1})$ and $(s_i,\ldots,s_k)$ of $s$, respectively. \hfill$\bullet$

In model predictive control, an input is typically obtained from the solution of a finite-horizon optimal control problem. The finite prediction horizon is some positive integer $N$ and the associated cost function is composed of a running cost $f_0\colon\mathbb{X}\times\mathbb{U}\to\mathbb{R}$ and a terminal cost $\phi\colon\mathbb{X}\to\mathbb{R}$. The cost function $J$ is then given by
\begin{equation}\label{eq:costFunction}
J(x_0,\nu) \ := \ \mathop\text{\footnotesize$\sum$}\limits_{k=0}^{N-1}f_0(x_k,\nu_k) + \phi(x_N)
\end{equation}
for every $x_0\in\mathbb{X}$ and every $\nu\in\mathbb{U}^{[\vphantom{]}0\colon\!\!N\vphantom{(})}$, where the states $x_k$ originate from a \emph{prediction model}
\begin{equation}\label{eq:predictionModel}
x_{k+1} \ = \ \hat{A} x_k + \hat{B} \nu_k.
\end{equation}
If the maps $A$ and $B$ in the actual control system \cref{eq:unknownSystem} are known, then $\hat{A}=A$ and $\hat{B}=B$ is the natural choice in \cref{eq:predictionModel}. If, however, $A$ and $B$ are unknown, then a suitable choice of $\hat{A}$ and $\hat{B}$ is, in general, not clear. Later, in \Cref{sec:4}, when we deal with unknown systems, we will use measurements of input-state pairs of \cref{eq:unknownSystem} to obtain estimates $(\hat{A},\hat{B})$ of $(A,B)$ for the prediction model \cref{eq:predictionModel}.

In standard MPC, only the first value of a predicted optimal input sequence is applied to the actual control system~\cref{eq:unknownSystem}. In contrast to this, flexible-step MPC allows for a flexible number of implemented values of an optimal input sequence \cite{Furnsinn2024}. The natural upper bound for the number of implement input values is the prediction horizon $N$ in \cref{eq:costFunction}. In \Cref{sec:3,sec:4}, we will present flexible-step algorithms for both known and unknown linear systems. The proposed algorithms will be composed of the following three components:
\begin{enumerate}[label=(C\arabic*),leftmargin=0.8cm]
	\item\label{item:C1} positive real numbers $\sigma_1,\ldots,\sigma_N$ such that
	\begin{equation}\label{eq:weights}
	\sigma_1 + \cdots + \sigma_N \ \geq \ 1;
	\end{equation}
	\item\label{item:C2} a positive real number $\alpha$ such that $0<\alpha<1$;
	\item\label{item:C3} a real-valued function $V$ on the state space $\mathbb{X}$ for which there exist positive real numbers $\underline{c}$, $\overline{c}$, and $p$ such that
	\begin{equation}\label{eq:Michel:eq:3.83}
	\underline{c}\,|x|^p \ \leq \ V(x) \ \leq \ \overline{c}\,|x|^p
	\end{equation}
	for every $x\in\mathbb{X}$, where $|\cdot|$ is a norm on $\mathbb{X}$.
\end{enumerate}
Following the terminology in \cite{Furnsinn2024}, the above function $V$ is said to be a \emph{generalized discrete-time control Lyapunov function} (GDCLF) for \cref{eq:predictionModel} if for every $x_0\in\mathbb{X}$, there exists $\nu\in\mathbb{U}^{[\vphantom{]}0\colon\!\!N\vphantom{(})}$ such that the so-called \emph{average descent condition}
\begin{equation}\label{eq:adc}
\mathop\text{\footnotesize$\sum$}\limits_{k=1}^N\sigma_k\,V(x_k) \ \leq \ (1-\alpha)\,V(x_0)
\end{equation}
is satisfied, where the states $x_1,\ldots,x_N\in\mathbb{X}$ are given by \cref{eq:predictionModel}. In \cref{eq:adc}, the positive real numbers $\sigma_1,\ldots,\sigma_N$ serve as \emph{weights} for the values of $V$ at the predicts states, and the positive real number $\alpha$ may be viewed as a \emph{decay constant}, which quantifies the decay of the weighted sum with respect to the initial value. If $V$ is a GDCLF for \cref{eq:predictionModel}, then one can show that $V$ is a non-monotonic Lyapunov function (see proof of \Cref{thm:Michel:Theorem3.5.7}), which in turn implies that \cref{eq:predictionModel} is asymptotically stable. We make the following simple but important observation for the purpose of flexible-step MPC.
\begin{remark}\label{rmk:flexibleStepExistence}
If the average descent condition~\cref{eq:adc} is satisfied for predicted states $x_1,\ldots,x_N$, then we can choose a ``flexible step'' number $\ell\in[1\colon\!\!N]$ such that the descent condition $V(x_\ell)-V(x_0)\leq-\alpha\,V(x_0)$ holds. The choice of a ``flexible step'' number will appear in the proposed flexible-step MPC \Cref{algo:KLTI,algo:ULTI} in \Cref{sec:3,sec:4}, respectively.\hfill$\bullet$
\end{remark}
Fix arbitrary components $\sigma_1,\ldots,\sigma_N$, $\alpha$, $V$ as in \ref{item:C1}-\ref{item:C3}. Then, for every $\hat{A}\in\mathcal{L}(\mathbb{X},\mathbb{X})$, every $\hat{B}\in\mathcal{L}(\mathbb{U},\mathbb{X})$, and every $\bar{x}\in\mathbb{X}$,  we consider the optimal control problem
\begin{align*}
\mathrm{OCP}(\hat{A},\hat{B},\bar{x})\colon \quad \mathop{\text{minimize}}\limits_{(x,\nu)} \ & \mathop\text{\footnotesize$\sum$}\limits_{k=0}^{N-1}f_0(x_k,\nu_k) + \phi(x_N) \\
\text{subject to}\, \ & (x,\nu)\in\Omega(\hat{A},\hat{B},\bar{x}),
\end{align*}
where minimization takes places over the set
\begin{align*}
\Omega(\hat{A},\hat{B},\bar{x}) \ := \ \Big\{ \ & (x,\nu)\in\mathbb{X}^{[0\colon\!\!N]}\times\mathbb{U}^{[\vphantom{]}0\colon\!\!N\vphantom{(})} \ \ \Big| \ \ x_0 = \bar{x}, \\
& \quad\ \forall{k}\in[\vphantom{]}0\colon\!\!N\vphantom{(})\colon x_{k+1} = \hat{A} x_k + \hat{B} \nu_k, \\
& \quad\ \mathop\text{\footnotesize$\sum$}\limits_{k=1}^N\sigma_k\,V(x_k) \ \leq \ (1-\alpha)\,V(x_0) \ \Big\}.
\end{align*}
The above optimal control problem $\mathrm{OCP}(\hat{A},\hat{B},\bar{x})$ is said to be \emph{feasible} if the set $\Omega(\hat{A},\hat{B},\bar{x})$ is not empty. The constraints in the feasible set consist of the dynamics of the prediction model \cref{eq:predictionModel} and the average descent condition \cref{eq:adc}. Clearly, $V$ is a GDCLF for \cref{eq:predictionModel} if and only if, for every $\bar{x}\in\mathbb{X}$, the optimal control problem $\mathrm{OCP}(\hat{A},\hat{B},\bar{x})$ is feasible.

\emph{From now on, we make the assumption that the optimal control problem $\mathrm{OCP}(\hat{A},\hat{B},\bar{x})$ has a solution whenever the set $\Omega(\hat{A},\hat{B},\bar{x})$ over which is optimized is not empty.}

The above assumption is satisfied under mild assumptions on the cost function; for example, if $f_0$ and $\phi$ are continuous, positive definite, and radially unbounded.
\begin{remark}\label{rmk:OptimizationStabilization}
Many MPC methods rely on a suitable design of the terminal costs (and/or additional terminal constraints) in order to achieve asymptotic stability. While the running cost $f_0$ is often part of the task specification, the choice of the terminal cost $\phi$ typically requires some prior knowledge about the to-be-controlled system. We will see that the proposed flexible-step approach does not involve any assumption about the choice of the cost functions $f_0$ and $\phi$, even in the case of an unknown control system. The problem of stabilization amounts to the problem of choosing a GDCLF for the unknown system. Since the GDCLF only appears in the constraints, the problem of optimization is, to some extent, decoupled from the problem of stabilization. Very simple functions, such as norm functions or quadratic functions, can be used as GDCLFs for entire classes of asymptotically stabilizable systems. \hfill$\bullet$
\end{remark}

\section{Flexible-step MPC for known LTI systems}\label{sec:3}
Suppose that the components $\sigma_1,\ldots,\sigma_N$, $\alpha$, $V$ are selected as in \ref{item:C1}-\ref{item:C3} of \Cref{sec:2}. In this section, we assume that the maps $A$ and $B$ in \cref{eq:unknownSystem} are known and therefore we can use them in the prediction model \cref{eq:predictionModel}. We will see below that the following assumption is sufficient to guarantee global exponential stability of \cref{eq:unknownSystem} under flexible-step MPC.
\begin{assumption}\label{ass:feasibility}
Suppose that, for every $\bar{x}\in\mathbb{X}$, the feasible set $\Omega(A,B,\bar{x})$ is not empty.
\end{assumption}
\begin{remark}[Choice of Parameters]\label{rmk:ChoiceOfParameters}
Notice that if control system \cref{eq:unknownSystem} is stabilizable, then it is possible to satisfy \Cref{ass:feasibility} through a suitable choice of the prediction horizon $N$ and the flexible-step components $\sigma_1,\ldots,\sigma_N$, $\alpha$, and $V$. To see this, recall that stabilizability of an LTI system is equivalent to the existence of some linear state feedback for which the closed-loop system is globally exponentially stable. For such a feedback law there exist sufficiently large $N\in\mathbb{N}$ and sufficiently small $\alpha\in(0,1)$ such that, for every solution $x\colon\mathbb{N}\to\mathbb{X}$ of the closed-loop system and every time $t\in\mathbb{N}$, the descent condition $|x(t+N)|\leq\frac{1-\alpha}{2}|x(t)|$ holds. Consequently, to satisfy \Cref{ass:feasibility}, we can select $V:=|\cdot|$, choose $\sigma_1,\ldots,\sigma_{N-1}>0$ sufficiently small, and set $\sigma_N:=1-\sum_{k=1}^{N-1}\sigma_k$. In particular (as a general guideline for the choice of parameters), if \cref{eq:unknownSystem} is stabilizable but unknown (as in \Cref{sec:4}), it is advisable to choose $N$ large, $\alpha,\sigma_1,\ldots,\sigma_{N-1}$ small, $\sigma_N:=1-\sum_{k=1}^{N-1}\sigma_k$, and $V:=|\cdot|$.

A direct and explicit choice of parameters is possible whenever the system \cref{eq:unknownSystem} is \emph{controllable}, that is, if the \emph{controllability matrix}
\begin{equation*}
\mathcal{C} \ := \ \begin{bmatrix} B & AB & \cdots & A^{n-1}B \end{bmatrix} \ \in \ \mathcal{L}(\mathbb{U}^n,\mathbb{X})
\end{equation*}
has maximum rank $n:=\dim\mathbb{X}$. In this case, for every state $x\in\mathbb{X}$, there exists an input sequence that steers the state of \cref{eq:unknownSystem} in at most $n$ steps from $x$ to the origin. One can compute such an input sequence by means of the pseudoinverse of $\mathcal{C}$. Suppose that a lower bound $r_{\mathcal{C}}>0$ of the smallest singular value of $\mathcal{C}$ is known and suppose that upper bounds $R_A>0$ and $R_B>0$ of the largest singular values of $A$ and $B$ are known. Then, one can choose $V:=|\cdot|$ and arbitrary $N\geq{n}$ and $\alpha\in(0,1)$. Based on the knowledge of the bounds $r_{\mathcal{C}}$, $R_A$, $R_B$, it is now straightforward to select weights as in \ref{item:C1} for which the average descent condition \cref{eq:adc} holds.
\hfill$\bullet$
\end{remark}
If $A$, $B$ are known and if \Cref{ass:feasibility} is satisfied, then the following algorithm is well-defined.
\begin{algorithm}[flexible-step MPC for known systems]\label{algo:KLTI} $ $
\begin{enumerate}[nolistsep,label=\arabic*:,ref=\arabic*,leftmargin=0.5cm]\setcounter{enumi}{-1}
	\item\label{algo:KLTI:0} Measure the initial state $x(0)$ of \cref{eq:unknownSystem} at initial time $t:=0$.
	\item\label{algo:KLTI:1} Compute a solution $(x,\nu)$ of $\mathrm{OCP}(A,B,x(t))$.
	\item\label{algo:KLTI:2} Choose a ``flexible step'' number $\ell\in[1\colon\!\!N]$ such that the descent condition $V(x_\ell)-V(x_0)\leq-\alpha\,V(x_0)$ holds.
	\item\label{algo:KLTI:3}
		\begin{enumerate}[nolistsep,label=\arabic{enumi}.\arabic*:,ref=\arabic{enumi}.\arabic*,leftmargin=1.0cm]\setcounter{enumii}{-1}
            \item\label{algo:KLTI:3.0} Set $k:=0$.
			\item\label{algo:KLTI:3.1} Apply the input $u(t):=\nu_k$ to \cref{eq:unknownSystem}.
			\item\label{algo:KLTI:3.2} Increment $k:=k+1$ and $t:=t+1$.
			\item\label{algo:KLTI:3.4} If $k<\ell$, go to \ref{algo:KLTI:3.1}, otherwise continue.
		\end{enumerate}
		\item\label{algo:KLTI:4} Measure the state $x(t)$ of \cref{eq:unknownSystem} and go to \ref{algo:KLTI:1}.
\end{enumerate}
\end{algorithm}
\begin{theorem}\label{thm:Michel:Theorem3.5.7}
Suppose \Cref{ass:feasibility} is satisfied. Then, system \cref{eq:unknownSystem} under \Cref{algo:KLTI} is globally exponentially stable.
\end{theorem}
\begin{proof}
We verify that $V$ is a non-monotonic Lyapunov function for the closed-loop system. Recall that we assume the existence of positive real numbers $\underline{c}$, $\overline{c}$, and $p$ such that \cref{eq:Michel:eq:3.83} is satisfied for every $x\in\mathbb{X}$. For the rest of the proof, fix an arbitrary solution $x\colon\mathbb{N}\to\mathbb{X}$ of system \cref{eq:unknownSystem} under \Cref{algo:KLTI}. Let $0=\tau_0<\tau_1<\cdots$ denote all time instants at which line~\ref{algo:KLTI:1} of \Cref{algo:KLTI} is executed. The choice of a ``flexible step'' number in line \ref{algo:KLTI:2} of \Cref{algo:KLTI} implies that $1\leq\tau_{\kappa+1}-\tau_\kappa\leq{N}$ and that the descent condition
\begin{equation*}
V(x(\tau_{\kappa+1})) - V(x(\tau_\kappa)) \ \leq \ - \alpha \cdot V(x(\tau_\kappa))
\end{equation*}
is satisfied for every $\kappa\in\mathbb{N}$. The decay of $V$ along the trajectory $x$ may be non-monotonic. However, since the predicted trajectory in line~\ref{algo:KLTI:1} of \Cref{algo:KLTI} satisfies the average descent condition \cref{eq:adc}, we also know that, for every $\kappa\in\mathbb{N}$ and every intermediate time instant $t\in(\tau_\kappa\colon\!\!\tau_{\kappa+1})$, the boundedness property
\begin{equation*}
V(x(t)) \ \leq \ \tfrac{1-\alpha}{\sigma_{\text{min}}} \cdot V(x(\tau_\kappa))
\end{equation*}
holds, where $\sigma_{\text{min}}>0$ is the minimum of the positive weights $\sigma_1,\ldots,\sigma_N$. Because of the above properties of $V$ along solutions of the closed-loop system, standard results for non-monotonic Lyapunov functions, such as Theorem~3.5.7 in \cite{MichelBook}, imply global exponential stability.
\end{proof}

\section{Flexible-step MPC for unknown LTI systems}\label{sec:4}
Suppose that the components $\sigma_1,\ldots,\sigma_N$, $\alpha$, $V$ are selected as in \ref{item:C1}-\ref{item:C3} of \Cref{sec:2}. In this section, the maps $A$ and $B$ in \cref{eq:unknownSystem} are not assumed to be known. That is, we are dealing with an unknown LTI system. We assume, however, that the inputs and states of \cref{eq:unknownSystem} can be measured at any time. This assumption will allow us to collect measured data in form of \emph{data lists}, which in turn are then used to obtain estimates $(\hat{A},\hat{B})$ of the unknown pair $(A,B)$.

\subsection{Estimates of the unknown system}\label{sec:4.1}
For the moment, fix an arbitrary time instant $\tau\in\mathbb{N}$. By a \emph{data list at time $\tau$}, we mean tuple $\mathcal{D}_\tau$ of the form
\begin{equation}\label{eq:dataList1}
\mathcal{D}_\tau \ = \ \big(x(0),(u(0),x(1)),\ldots,(u(\tau-1),x(\tau))\big),
\end{equation}
where the entries $x(t)\in\mathbb{X}$ and $u(t)\in\mathbb{U}$ satisfy \cref{eq:unknownSystem} for every $t\in[\vphantom{]}0\colon\!\!\tau\vphantom{(})$. After applying an input $u(\tau)\in\mathbb{U}$ to \cref{eq:unknownSystem} and measuring the next state $x(\tau+1)\in\mathbb{X}$, we can \emph{append} the input-state pair $(u(\tau),x(\tau+1))$ to $\mathcal{D}_\tau$ and obtain the new data list
\begin{equation}\label{eq:dataList2}
\mathcal{D}_{\tau+1} \ = \ \big(x(0),(u(0),x(1)),\ldots,(u(\tau),x(\tau+1))\big)
\end{equation}
at time $\tau+1$. For every data list $\mathcal{D}_\tau$ of the form \cref{eq:dataList1}, define $\mathcal{E}(\mathcal{D}_\tau)$ to be the set of all pairs $(\hat{A},\hat{B})$ with $\hat{A}\in\mathcal{L}(\mathbb{X},\mathbb{X})$ and $\hat{B}\in\mathcal{L}(\mathbb{U},\mathbb{X})$ such that
\begin{equation}\label{eq:trajectoryCondition}
x(t+1) \ = \ \hat{A} x(t) + \hat{B} u(t)
\end{equation}
for every $t\in[\vphantom{]}0\colon\!\!\tau\vphantom{(})$. The set $\mathcal{E}(\mathcal{D}_\tau)$ contains all potential estimates of the unknown pair $(A,B)$ based on the knowledge provided by the data list $\mathcal{D}_\tau$. Since $\mathcal{E}(\mathcal{D}_{\tau+1})$ is a subset of $\mathcal{E}(\mathcal{D}_\tau)$, the following definition makes sense.
\begin{definition}\label{def:estimator}
By an \emph{estimator of $(A,B)$}, we mean a map $E$ that assigns to every data set $\mathcal{D}_{\tau+1}$ of the form \cref{eq:dataList2} a map
\begin{equation}\label{eq:estimator1}
E(\mathcal{D}_{\tau+1})\colon \mathcal{E}(\mathcal{D}_{\tau}) \to \mathcal{E}(\mathcal{D}_{\tau+1})
\end{equation}
such that
\begin{equation}\label{eq:estimator2}
E(\mathcal{D}_{\tau+1})(\hat{A},\hat{B}) \ = \ (\hat{A},\hat{B})
\end{equation}
for every $(\hat{A},\hat{B})\in\mathcal{E}(\mathcal{D}_{\tau+1})$, where $\mathcal{D}_{\tau}$ is the truncated data list of $\mathcal{D}_{\tau+1}$ of the form \cref{eq:dataList1}. \hfill$\bullet$
\end{definition}
\begin{remark}
Suppose that $E$ is an estimator of $(A,B)$ and that $\mathcal{D}_{\tau}$ and $\mathcal{D}_{\tau+1}$ are two consecutive data lists of the form \cref{eq:dataList1} and \cref{eq:dataList2}. Then, the map $E(\mathcal{D}_{\tau+1})$ in \cref{eq:estimator1} assigns to any given estimate of $(A,B)$ at time $\tau$ the next estimate of $(A,B)$ at time $\tau+1$, which takes the gain of information about the unknown system into account. Equation \cref{eq:estimator2} means that the estimate remains unchanged if the previous estimate can also be used as the next one. \hfill$\bullet$
\end{remark}
\begin{example}\label{exm:estimator}
A simple example of an estimator $E$ as in \Cref{def:estimator} is the following \emph{least norm} (or \emph{least squares}) estimator. Let $|\cdot|_{\mathbb{X},\mathbb{X}}$ and $|\cdot|_{\mathbb{U},\mathbb{X}}$ be norms on $\mathcal{L}(\mathbb{X},\mathbb{X})$ and $\mathcal{L}(\mathbb{U},\mathbb{X})$, respectively, which are assumed to be induced by inner products. Suppose that $\mathcal{D}_{\tau}$ and $\mathcal{D}_{\tau+1}$ are two consecutive data lists of the form \cref{eq:dataList1} and \cref{eq:dataList2}. Then, we define the map $E(\mathcal{D}_{\tau+1})$ in \cref{eq:estimator1} as follows. For every $(\hat{A},\hat{B})\in\mathcal{E}(\mathcal{D}_{\tau+1})$, let \cref{eq:estimator2} be the defining rule of $E(\mathcal{D}_{\tau+1})$. For $(\hat{A},\hat{B})\notin\mathcal{E}(\mathcal{D}_{\tau+1})$, define $E(\mathcal{D}_{\tau+1})(\hat{A},\hat{B})$ to be the unique minimizer of the function
\begin{equation}\label{eq:LS}
\mathcal{E}(\mathcal{D}_{\tau+1}) \to \mathbb{R}, \qquad (\tilde{A},\tilde{B}) \mapsto |\tilde{A}|_{\mathbb{X},\mathbb{X}}^2 + |\tilde{B}|_{\mathbb{U},\mathbb{X}}^2.
\end{equation}
Then $E$ is an estimator of $(A,B)$. \hfill$\bullet$
\end{example}

\subsection{Persistently exciting inputs}\label{sec:4.2}
Recall the following definition of persistently exciting inputs from \cite{Willems2005}.
\begin{definition}\label{def:PE}
For any positive integers $T$ and $d$, a finite se\-quence $v\in\mathbb{U}^{[\vphantom{]}0\colon\!\!T\vphantom{(})}$ is \emph{persistently exciting of order~$d$} if the subsequences $v_{[\vphantom{]}0\colon\!\!d\vphantom{(})}$, $\ldots$, $v_{[\vphantom{]}T-d\colon\!\!T\vphantom{(})}$ of $v$ span $\mathbb{U}^d$. \hfill$\bullet$
\end{definition}
In the proposed flexible-step MPC approach for unknown linear systems (see \Cref{algo:ULTI} below), a persistently exciting input sequence of order $\dim\mathbb{X}+1$ will be used as a last resort in case of an infeasible optimal control problem to explore the unknown system. Such a sequence can be easily generated as follows.
\begin{example}\label{exm:PE}
Let $d:=\dim\mathbb{X}+1$ and $m:=\dim\mathbb{U}$. Choose a basis $\{b_1,\ldots,b_m\}$ of $\mathbb{U}$ and an integer $T\geq(m+1)d-1$. Define $v\in\mathbb{U}^{[\vphantom{]}0\colon\!\!T\vphantom{(})}$ by
\begin{align*}
v_k & := \left\{ \begin{tabular}{cl} $b_i$ & if $k=i\cdot{d}-1$ for some $i\in\{1,\ldots,m\}$, \\ $0$ & otherwise. \end{tabular} \right.
\end{align*}
Then, $v$ is persistently exciting of order $d$. As an alternative, for $T$ as above, one can also choose $v_0,\ldots,v_{T-1}$ as independent random vectors with Gaussian distribution. Then, \Cref{lem:PE} in \hyperlink{sec:appendix}{Appendix} guarantees that the so-defined random sequence $v$ of length $T$ is almost surely persistently exciting of order $d$. \hfill$\bullet$
\end{example}
In the next remark, we indicate the role of persistently exciting inputs in our flexible-step MPC approach. 
\begin{remark}\label{rmk:Yu}
Let $E$ be an estimator of $(A,B)$ as in \Cref{def:estimator}. Suppose that, at some time instant $\tau\in\mathbb{N}$, we have gathered a data list $\mathcal{D}_\tau$ of the form \cref{eq:dataList1}. Let $(\hat{A}_\tau,\hat{B}_\tau)\in\mathcal{E}(\mathcal{D}_\tau)$ be the current estimate of $(A,B)$ based on the information provided by $\mathcal{D}_\tau$. Choose a persistently exciting sequence $v\in\mathbb{U}^{[\vphantom{]}0\colon\!\!T\vphantom{(})}$ of order $\dim\mathbb{X}+1$. Then, inductively, at every time instant $t\in[\vphantom{]}\tau\colon\!\!\tau+T\vphantom{(})$,
\begin{enumerate}[nolistsep,label=\arabic*:,ref=\arabic*,leftmargin=0.5cm]
	\item apply the input $u(t):=v_{t-\tau}$ to \cref{eq:unknownSystem},
	\item measure the next state $x(t+1)$ of \cref{eq:unknownSystem},
	\item append the pair $(u(t),x(t+1))$ to $\mathcal{D}_t$ to obtain $\mathcal{D}_{t+1}$,
	\item update the estimate $(\hat{A}_{t+1},\hat{B}_{t+1}):=E(\mathcal{D}_{t+1})(\hat{A}_t,\hat{B}_t)$.
\end{enumerate}
By the construction, the estimate $(\hat{A}_{\tau+T},\hat{B}_{\tau+T})$ satisfies
\begin{equation}\label{eq:Yu1}
x(t+1) \ = \ \hat{A}_{\tau+T} x(t) + \hat{B}_{\tau+T} v_{t-\tau}
\end{equation}
for every $t\in[\vphantom{]}\tau\colon\!\!\tau+T\vphantom{(})$. Since $v$ is assumed to be persistently exciting of order $\dim\mathbb{X}+1$, we may utilize Theorem~1 of \cite{Yu2021}, which is a generalization of Willems' fundamental lemma. This result implies that, for every subsequently applied input sequence $u\colon[\vphantom{]}\tau+T\colon\!\!\infty\vphantom{(})\to\mathbb{U}$ and respective state sequence $x\colon[\vphantom{]}\tau+T\colon\!\!\infty\vphantom{(})\to\mathbb{X}$ of \cref{eq:unknownSystem}, the estimate $(\hat{A}_{\tau+T},\hat{B}_{\tau+T})$ satisfies
\begin{equation}\label{eq:Yu2}
x(t+1) \ = \ \hat{A}_{\tau+T} x(t) + \hat{B}_{\tau+T} u(t)
\end{equation}
at every time $t\geq\tau+T$. In other words, after an application of $v$, the resulting data list provides enough information to predict all future states of the unknown system \cref{eq:unknownSystem}. Since the estimator $E$ is assumed to satisfy \cref{eq:estimator2}, it follows that
\begin{equation}\label{eq:Yu3}
(\hat{A}_{\tau+T},\hat{B}_{\tau+T}) = (\hat{A}_{t+1},\hat{B}_{t+1}) := E(\mathcal{D}_{t+1})(\hat{A}_t,\hat{B}_t)
\end{equation} 
for every time $t\geq\tau+T$; that is, the estimates remain unchanged after time $\tau+T$. More importantly, for our objectives regarding flexible-step MPC, it follows that, for every state $x(t)$ of \cref{eq:unknownSystem} at some subsequent time instant $t\geq\tau+T$, the optimal control problems $\mathrm{OCP}(\hat{A}_t,\hat{B}_t,x(t))$ and $\mathrm{OCP}(A,B,x(t))$ are equivalent in the sense that their feasible sets
\begin{equation}\label{eq:Yu4}
\Omega(\hat{A}_t,\hat{B}_t,x(t)) \ = \ \Omega(A,B,x(t))
\end{equation}
coincide. Consequently, from time $\tau+T$ on, the action of \Cref{algo:KLTI} on system \cref{eq:unknownSystem} remains unchanged if the unknown pair $(A,B)$ in \Cref{algo:KLTI} is replaced by the known estimates. In particular, if \Cref{ass:feasibility} is satisfied, then, by \Cref{thm:Michel:Theorem3.5.7}, an application of \Cref{algo:KLTI} with the estimates of $(A,B)$ will lead to exponential convergence of the system state to the origin.

Of course, it would be also possible to solve the control problem by first applying an entire persistently exciting input sequence of order $\dim\mathbb{X}+1$ to explore \cref{eq:unknownSystem}, and then employing \Cref{algo:KLTI} for known linear systems with reliable estimates of $(A,B)$. Such an approach is taken, for example, in \cite{Berberich20202}. In practice, however, this may not be ideal, because an application of an entire persistently exciting input sequence can lead to an undesirable system behavior, or may be of high cost. For this reason, we will only use persistently exciting inputs if absolutely necessary. Namely, in case that the current prediction model leads to an infeasible optimal control problem. As soon as the optimal control problem becomes feasible again, we use flexible-step MPC for simultaneous exploration and exploitation. Notably, the simulation results in \Cref{sec:5} confirm that the flexible-step approach can lead to a better system behavior than a naive application of persistently exciting inputs for exploration.\hfill$\bullet$
\end{remark}

\subsection{Proposed algorithm and main result}\label{sec:4.3}
Regardless of whether \Cref{ass:feasibility} is satisfied, the following algorithm is always well-defined.
\begin{algorithm}[flexible-step MPC for unknown systems]\label{algo:ULTI}$ $
\begin{enumerate}[nolistsep,label=\arabic*:,ref=\arabic*,leftmargin=0.5cm]\setcounter{enumi}{-1}
	\item\label{algo:ULTI:0} Choose an estimator $E$ of $(A,B)$ as in \Cref{def:estimator}. Choose an arbitrary initial estimate $(\hat{A}_0,\hat{B}_0)$ of $(A,B)$. Measure the initial state $x(0)$ of \cref{eq:unknownSystem} at initial time $t:=0$. Define the initial data list $\mathcal{D}_0:=(x(0))$. 
	\item\label{algo:ULTI:1} If $\Omega(\hat{A}_t,\hat{B}_t,x(t))$ is not empty, go to \ref{algo:ULTI:2}, otherwise to \ref{algo:ULTI:3}.
	\item\label{algo:ULTI:2} \begin{enumerate}[nolistsep,label=\arabic{enumi}.\arabic*:,ref=\arabic{enumi}.\arabic*,leftmargin=0.7cm]
		\item\label{algo:ULTI:2.1} Compute a solution $(x,\nu)$ of $\mathrm{OCP}(\hat{A}_t,\hat{B}_t,x(t))$.
		\item\label{algo:ULTI:2.2} Choose a ``flexible step'' number $\ell\in[1\colon\!\!N]$ such that the descent condition $V(x_\ell)-V(x_0)\leq-\alpha\,V(x_0)$ holds.
		\item\label{algo:ULTI:2.3}
		\begin{enumerate}[nolistsep,label=\arabic{enumi}.\arabic{enumii}.\arabic*:,ref=\arabic{enumi}.\arabic{enumii}.\arabic*,leftmargin=1.0cm]\setcounter{enumiii}{-1}
            \item\label{algo:ULTI:2.3.0} Set $k:=0$.
			\item\label{algo:ULTI:2.3.1} Apply the input $u(t):=\nu_k$ to \cref{eq:unknownSystem}.
			\item\label{algo:ULTI:2.3.2} Increment $k:=k+1$ and $t:=t+1$.
			\item\label{algo:ULTI:2.3.3} Measure the state $x(t)$ of \cref{eq:unknownSystem}.
			\item\label{algo:ULTI:2.3.4} Append $(u(t-1),x(t))$ to $\mathcal{D}_{t-1}$ to obtain $\mathcal{D}_t$.
            \item\label{algo:ULTI:2.3.5} Set $(\hat{A}_t,\hat{B}_t):=E(\mathcal{D}_t)(\hat{A}_{t-1},\hat{B}_{t-1})$.
			\item\label{algo:ULTI:2.3.6} If $k<\ell$ go to \ref{algo:ULTI:2.3.1}, otherwise go to \ref{algo:ULTI:1}.
		\end{enumerate}
	\end{enumerate}
	\item\label{algo:ULTI:3}
	\begin{enumerate}[nolistsep,label=\arabic{enumi}.\arabic*:,ref=\arabic{enumi}.\arabic*,leftmargin=0.7cm]\setcounter{enumii}{-1}
		\item\label{algo:ULTI:3.0} Set $k:=0$ and $\tau:=t$. Choose a persistently exciting sequence $v\in\mathbb{U}^{[\vphantom{]}0\colon\!\!T\vphantom{(})}$ of order $\dim\mathbb{X}+1$.
		\item\label{algo:ULTI:3.1} Apply the input $u(t):=v_k$ to \cref{eq:unknownSystem}.
		\item\label{algo:ULTI:3.2} Increment $k:=k+1$ and $t:=t+1$.
		\item\label{algo:ULTI:3.3} Measure the state $x(t)$ of \cref{eq:unknownSystem}.
		\item\label{algo:ULTI:3.4} Append $(u(t-1),x(t))$ to $\mathcal{D}_{t-1}$ to obtain $\mathcal{D}_t$.
		\item\label{algo:ULTI:3.5} Set $(\hat{A}_t,\hat{B}_t):=E(\mathcal{D}_t)(\hat{A}_{t-1},\hat{B}_{t-1})$.
		\item\label{algo:ULTI:3.6} If $k<T$ and if $(\hat{A}_t,\hat{B}_t)=(\hat{A}_\tau,\hat{B}_\tau)$, go to \ref{algo:ULTI:3.1}, otherwise go to \ref{algo:ULTI:1}.
	\end{enumerate}
\end{enumerate}
\end{algorithm}
We can prove the following global convergence property of the unknown system \cref{eq:unknownSystem} under \Cref{algo:ULTI}.
\begin{theorem}\label{thm:mainResult}
Suppose that \Cref{ass:feasibility} is satisfied. Then, for every initial point $x(0)$, the state $x(t)$ of system \cref{eq:unknownSystem} under \Cref{algo:ULTI} converges to the origin as $t\to\infty$.
\end{theorem}
\begin{proof}
Fix an arbitrary initial state $x(0)\in\mathbb{X}$, an arbitrary initial estimate $(\hat{A}_0,\hat{B}_0)$ of $(A,B)$, and set $\mathcal{D}_0:=(x(0))$. Then, for every $t\in\mathbb{N}$, \Cref{algo:ULTI} generates well-defined $\hat{A}_t\in\mathcal{L}(\mathbb{X},\mathbb{X})$, $\hat{B}_t\in\mathcal{L}(\mathbb{U},\mathbb{X})$, $x(t)\in\mathbb{X}$, $u(t)\in\mathbb{U}$, and $\mathcal{D}_t$ of the form \cref{eq:dataList1}.

First, we treat the case in which there exists a time instant $\tau\in\mathbb{N}$ at which line \ref{algo:ULTI:3.0} of \Cref{algo:ULTI} is executed, and then lines \ref{algo:ULTI:3.1}-\ref{algo:ULTI:3.6} are repeatedly executed for the entire length $T$ of the persistently exciting sequence $v$ order $\dim\mathbb{X}+1$. In this case, we know from \Cref{rmk:Yu} that \cref{eq:Yu4} holds at every time $t\geq\tau+T$. Since \Cref{ass:feasibility} is satisfied, we conclude that $\Omega(\hat{A}_t,\hat{B}_t,x(t))$ is not empty for every $t\geq\tau+T$. It follows that only the lines in part~\ref{algo:ULTI:2} of \Cref{algo:ULTI} are executed from time $\tau+T$ on. Again, since \cref{eq:Yu4} holds for every $t\geq\tau+T$, part~\ref{algo:ULTI:2} of \Cref{algo:ULTI} has the same effect on \cref{eq:unknownSystem} as \Cref{algo:KLTI}. Therefore, \Cref{thm:Michel:Theorem3.5.7} implies the claim.

Because of the considerations in the previous paragraph, we may assume from now on that, for every time instant $\tau\in\mathbb{N}$ at which line \ref{algo:ULTI:3.0} of \Cref{algo:ULTI} is executed, there exists $t\in[\vphantom{]}\tau\colon\!\!\tau+T\vphantom{(})$ such that the loop consisting of lines \ref{algo:ULTI:3.1}-\ref{algo:ULTI:3.6} is terminated at time $t$ due to $(\hat{A}_t,\hat{B}_t)\neq(\hat{A}_\tau,\hat{B}_\tau)$. Notice that $\mathcal{E}(\mathcal{D}_0)\supseteq\mathcal{E}(\mathcal{D}_1)\supseteq\cdots$ is a non-increasing sequence of affine vector spaces. Since the employed estimator $E$ has the property in equation \cref{eq:estimator2}, it follows that, for every $t\in\mathbb{N}$, the following implication holds: If $(\hat{A}_{t+1},\hat{B}_{t+1})$ deviates from $(\hat{A}_t,\hat{B}_t)$, then $\mathcal{E}(\mathcal{D}_{t+1})$ has lower dimension than $\mathcal{E}(\mathcal{D}_t)$. This in turn implies that there exists some sufficiently large $\tau\in\mathbb{N}$ such that $(\hat{A}_t,\hat{B}_t)=(\hat{A}_\tau,\hat{B}_\tau)$ at every time $t\geq\tau$. Because of the above assumption about part \ref{algo:ULTI:3} of \Cref{algo:ULTI}, it follows that only the lines in part~\ref{algo:ULTI:2} are executed from time $\tau$ on. Moreover, since $(\hat{A}_t,\hat{B}_t)=(\hat{A}_\tau,\hat{B}_\tau)$, we have
\[
x(t+1) \ = \ \hat{A}_\tau x(t) + \hat{B}_\tau u(t)
\]
at every time $t\geq\tau$, which means that the prediction model provides the actual states of the unknown system. Thus, from time $\tau$ on, part~\ref{algo:ULTI:2} of \Cref{algo:ULTI} has the same effect on \cref{eq:unknownSystem} as \Cref{algo:KLTI}. Again, \Cref{thm:Michel:Theorem3.5.7} implies the claim.
\end{proof}

\section{Simulation results}\label{sec:5}
In this section, we test the proposed \Cref{algo:ULTI} for unknown LTI systems through numerical simulations. To this end, we consider a class of systems, which, according to \cite{Zeng2023}, is hard to learn to stabilize. This class consists of single-input systems with state space $\mathbb{X}=\mathbb{R}^n$ and input space $\mathbb{U}=\mathbb{R}$. Following \cite{Zeng2023}, the unknown system matrices $A$ and $B$ in \cref{eq:unknownSystem} are of the form
\begin{equation}\label{eq:sim-sys}
A \ := \ \begin{bmatrix}
   r   &   v   &   0   &\cdots &   0   \\
   0   &   0   &   v   &       &   0   \\
\vdots &\vdots &       &\ddots &       \\
   0   &   0   &   0   &       &   v   \\
   0   &   0   &   0   &\cdots &   0
\end{bmatrix} \quad \text{and} \quad B \ := \ \begin{bmatrix}
b \\ 0 \\ \vdots \\ 0 \\ v
\end{bmatrix}
\end{equation}
with parameters $r>1$, $0<v<(r-1)/2$ and $b\geq0$. Here, we choose $n:=7$, $r:=2$, $v:=0.8$, and $b:=0.7$.

The running cost $f_0$ and the terminal cost $\phi$ in the optimal control problem (cf.~\Cref{sec:2}) are given by $f_0(x,u):=u^2$ and $\phi(x):=0$. The prediction horizon is $N:=10$. Clearly, for this choice of $f_0$ and $\phi$, it is unlikely that a standard MPC approach would lead to an asymptotically stable closed-loop system. However, as we will see below, our flexible-step approach with a GDCLF in the constraints can stabilize the system since the tasks of optimization and stabilization are decoupled (cf.~\Cref{rmk:OptimizationStabilization}). A candidate $V$ for a GDCLF is given by $V(x):=|x|$, where $|\cdot|$ denotes the Euclidean norm. The weights as well as the decay constant for the average descent condition \cref{eq:adc} are chosen to be $\sigma_k:=0.001$ for every $k\in[0\colon\!\!N-1]$ and $\sigma_N:=0.991$ as well as $\alpha:=0.001$. Then, the components $\sigma_1,\ldots,\sigma_N$, $\alpha$, $V$ are selected as in \ref{item:C1}-\ref{item:C3} of \Cref{sec:2}.

For the estimator $E$ of $(A,B)$, which is used in lines \ref{algo:ULTI:2.3.5} and \ref{algo:ULTI:3.5} of \Cref{algo:ULTI}, we choose the least norm (or least squares) estimator from \Cref{exm:estimator} with respect to the Frobenius norm. In this case, the minimizer $(\hat{A}_{\tau+1},\hat{B}_{\tau+1})$ of \cref{eq:LS} can be computed in matrix form by the formula
\[
\begin{bmatrix} \hat{A}_{\tau+1} \!\!&\!\! \hat{B}_{\tau+1} \end{bmatrix} \, = \, \begin{bmatrix} x(1) \!&\! \cdots \!&\! x(\tau+1) \end{bmatrix}\begin{bmatrix} x(0) \!&\! \cdots \!&\! x(\tau) \\ u(0) \!&\! \cdots \!&\! u(\tau) \end{bmatrix}^+\!\!\!,
\]
where $(\cdot)^+$ denotes the Moore--Penrose inverse (or pseudoinverse). The exploration inputs $v_k\in\mathbb{R}$ in line \ref{algo:ULTI:3.1} of \Cref{algo:ULTI} are generated by independent and normally distributed random variables with zero mean and a variance of~$0.01$. Then, it follows from \Cref{exm:PE} and \Cref{lem:PE} that a sequence of at least $15$ of these random inputs is almost surely persistently exciting of order $8$.
\begin{figure}
    \centering
    \includegraphics[width=\linewidth]{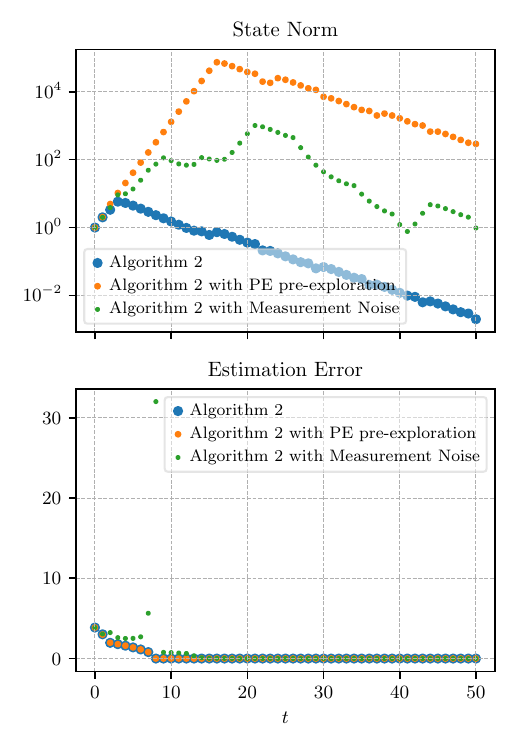}
    \caption{Simulation results for an unknown LTI system under flexible-step MPC. Upper plot: Euclidean norm of the system state as a function of time. Lower plot: Deviation of the prediction model from the unknown system. The blue dots represent the results generated by \Cref{algo:ULTI}. For comparison, the orange dots represent the results for the same choice of parameters, but with an additional exploration of the unknown system by persistently exciting (PE) inputs prior to the application of \Cref{algo:ULTI}. The pre-exploration by PE inputs lasts from the initial time $t=0$ until $t=15$. A direct application of \Cref{algo:ULTI} without pre-exploration (blue dots) only involves one PE input value at time $t=1$. For further comparison, the green dots represent the results generated by \Cref{algo:ULTI} for the same choice of parameters but under the additional influence of measurement noise. Again, only one PE input at time $t=1$ is applied.}
    \label{fig:simulation_comparison}
\end{figure}

For the initialization of \Cref{algo:ULTI}, we choose the initial system state $x(0) := [1,0,\ldots,0]^\top$ and the initial estimates $\hat{A}_0 := 0_{7\times 7}$ and $\hat{B}_0 := 0_7$ of $A$ and $B$. The``flexible step'' number $\ell$ in line \ref{algo:ULTI:2.2} of \Cref{algo:ULTI} is always selected as the one with the smallest value of $V$ along the predicted trajectory. The results generated by \Cref{algo:ULTI} are depicted in \Cref{fig:simulation_comparison} by blue dots. The upper plot in \Cref{fig:simulation_comparison} shows the Euclidean norm of the system state $x(t)$ as a function of the time $t$. It can be observed that exponential convergence sets in after a short transient of $4$ time steps. The lower plot in \Cref{fig:simulation_comparison} shows the deviation $|\hat{A}_t-A|_F + |\hat{B}_t-B|_F$ of the estimated pair $(\hat{A}_t,\hat{B}_t)$ in the prediction model from the unknown pair $(A,B)$ in the control system with respect to the Frobenius norm $|\cdot|_F$. It can be observed that full system identification is achieved after 9 time steps.

A beneficial feature of the proposed \Cref{algo:ULTI} is that persistently exciting (PE) inputs are only applied when absolutely necessary. As discussed in \Cref{rmk:Yu}, an application of PE inputs can be energy consuming and may result in an undesired system behavior. To visualize this fact, the orange dots in \Cref{fig:simulation_comparison} represent results for a naive PE approach. That is, prior to flexible-step MPC a full sequence of 15 PE inputs is fed into the system. After these 15 steps, the measurements of input-state pairs provide enough information to correctly predict all future states of the unknown system (cf.~\Cref{rmk:Yu}). The orange dots in the upper plot of \Cref{fig:simulation_comparison} indicate a significant increase of the state norm during the long exploration phase. An application of flexible-step MPC after the first 15 time steps then leads to exponential convergence of the system state to the origin. In contrast to the naive PE approach, the results generated by \Cref{algo:ULTI} (which are illustrated by blue dots in \Cref{fig:simulation_comparison}) involve only \emph{one} PE input value. This PE input is applied at time $t=1$, where part~\ref{algo:ULTI:3} of \Cref{algo:ULTI} is active. All other input values originate from the flexible-step routine in part~\ref{algo:ULTI:2} of \Cref{algo:ULTI}.

To evaluate the performance of \Cref{algo:ULTI} under more adverse conditions, we repeat the same simulation but now assume that the state measurements in lines \ref{algo:ULTI:0}, \ref{algo:ULTI:2.3}, and \ref{algo:ULTI:3.3} of \Cref{algo:ULTI} are corrupted by additive, uncorrelated Gaussian noise with zero mean and a standard deviation $0.05$. The corresponding results are depicted by green dots in \Cref{fig:simulation_comparison}. Compared to the noise-free case, the state norm increases for more time steps at the start of the simulation and the estimation error reaches a larger peak at $t=8$. Nevertheless, the results indicate that \Cref{algo:ULTI} has the ability to stabilize the system despite noisy measurements.

\section{Conclusions}\label{sec:6}
We presented a new flexible-step MPC approach that can be applied to unknown LTI systems without an unnecessarily long pre-exploration phase. The tasks of optimization and stabilization are decoupled through the introduction of a stabilizing \emph{generalized discrete-time control Lyapunov function} in the constraints. The employed cost function in the optimal control problem does not require a carefully designed terminal cost term. It can be freely chosen to enforce performance requirements and to include soft input and state constraints. We have seen in numerical simulations that the method also performs well in the presence of noise. A rigorous robustness analysis is left to future research. The flexible-step approach also bears some potential for extensions to unknown nonlinear systems. This includes nonlinear control systems whose state transition maps are \emph{unknown} linear combinations of \emph{known} nonlinear base functions.

\bibliographystyle{IEEEtran}
\bibliography{bibFile}

\begin{thebibliography}{10}
\providecommand{\url}[1]{#1}
\csname url@samestyle\endcsname
\providecommand{\newblock}{\relax}
\providecommand{\bibinfo}[2]{#2}
\providecommand{\BIBentrySTDinterwordspacing}{\spaceskip=0pt\relax}
\providecommand{\BIBentryALTinterwordstretchfactor}{4}
\providecommand{\BIBentryALTinterwordspacing}{\spaceskip=\fontdimen2\font plus
\BIBentryALTinterwordstretchfactor\fontdimen3\font minus
  \fontdimen4\font\relax}
\providecommand{\BIBforeignlanguage}[2]{{%
\expandafter\ifx\csname l@#1\endcsname\relax
\typeout{** WARNING: IEEEtran.bst: No hyphenation pattern has been}%
\typeout{** loaded for the language `#1'. Using the pattern for}%
\typeout{** the default language instead.}%
\else
\language=\csname l@#1\endcsname
\fi
#2}}
\providecommand{\BIBdecl}{\relax}
\BIBdecl

\bibitem{GruneBook}
L.~Gr\"une and J.~Pannek, \emph{{Nonlinear Model Predictive Control}},
  2nd~ed.\hskip 1em plus 0.5em minus 0.4em\relax Cham: Springer, 2017.

\bibitem{Hewing20202}
L.~Hewing, K.~Wabersich, M.~Menner, and M.~Zeilinger, ``{Learning-Based Model
  Predictive Control: Toward Safe Learning in Control},'' \emph{Annu. Rev.
  Control Robot. Auton. Syst.}, vol.~3, pp. 269--296, 2020.

\bibitem{Marafioti2014}
G.~Marafioti, R.~R. Bitmead, and M.~Hovd, ``{Persistently exciting model
  predictive control},'' \emph{Int. J. Adapt. Control Signal Process}, vol.~28,
  no.~6, pp. 536--552, 2014.

\bibitem{Heirung2017}
T.~A.~N. Heirung, B.~E. Ydstie, and B.~Foss, ``{Dual adaptive model predictive
  control},'' \emph{Automatica}, vol.~80, pp. 340--348, 2017.

\bibitem{Goncalves2016}
G.~A. Gon\c{c}alves and M.~Guay, ``{Robust discrete-time set-based adaptive
  predictive control for nonlinear systems},'' \emph{J. Process Control},
  vol.~39, pp. 111--122, 2016.

\bibitem{Aswani2013}
A.~Aswani, H.~Gonzalez, S.~S. Sastry, and C.~Tomlin, ``{Provably safe and
  robust learning-based model predictive control},'' \emph{Automatica},
  vol.~49, no.~5, pp. 1216--1226, 2013.

\bibitem{Tanaskovic2014}
M.~Tanaskovic, L.~Fagiano, R.~Smith, and M.~Morari, ``{Adaptive receding
  horizon control for constrained MIMO systems},'' \emph{Automatica}, vol.~50,
  no.~12, pp. 3019--3029, 2014.

\bibitem{Kim2008}
T.-H. Kim and T.~Sugie, ``{Adaptive receding horizon predictive control for
  constrained discrete-time linear systems with parameter uncertainties},''
  \emph{Int. J. Control}, vol.~81, no.~1, pp. 62--73, 2008.

\bibitem{WaardeBook}
H.~J. {van Waarde}, M.~K. Camlibel, and H.~L. Trentelman, \emph{{Data-Based
  Linear Systems and Control Theory}}.\hskip 1em plus 0.5em minus 0.4em\relax
  Seattle: Kindle Direct Publishing, 2025.

\bibitem{Coulson2019}
J.~Coulson, J.~Lygeros, and F.~D\"orfler, ``{Data-Enabled Predictive Control:
  In the Shallows of the DeePC},'' in \emph{{ECC}}, 2019, pp. 307--312.

\bibitem{Coulson20192}
------, ``Regularized and distributionally robust data-enabled predictive
  control,'' in \emph{{Proc. CDC}}, 2019, pp. 2696--2701.

\bibitem{Baros2022}
S.~Baros, C.-Y. Chang, G.~E. Col\'{o}n-Reyes, and A.~Bernstein, ``Online
  data-enabled predictive control,'' \emph{Automatica}, vol. 138, p. 109926,
  2022.

\bibitem{Berberich20202}
J.~Berberich, A.~Koch, M.~A. M\"uller, and F.~Allg\"ower, ``{Data-Driven Model
  Predictive Control With Stability and Robustness Guarantees},'' \emph{IEEE
  Trans. Automat. Control}, vol.~66, no.~4, pp. 1702--1717, 2021.

\bibitem{Berkenkamp2017}
F.~Berkenkamp, M.~Turchetta, A.~Schoellig, and A.~Krause, ``{Safe Model-based
  Reinforcement Learning with Stability Guarantees},'' in \emph{Proc. NeurIPS},
  2017, pp. 908--918.

\bibitem{Patan2015}
K.~Patan, ``{Neural Network-Based Model Predictive Control: Fault Tolerance and
  Stability},'' \emph{IEEE Trans. Control Syst. Technol.}, vol.~23, no.~3, pp.
  1147--1155, 2015.

\bibitem{Maiworm2018}
M.~Maiworm, D.~Limon, J.~{Maria Manzano}, and R.~Findeisen, ``{Stability of
  Gaussian Process Learning Based Output Feedback Model Predictive Control},''
  in \emph{Proc. NMPC}, 2018, pp. 455--461.

\bibitem{Yu20232}
J.~Yu, D.~Ho, and A.~Wierman, ``{Online Adversarial Stabilization of Unknown
  Linear Time-Varying Systems},'' in \emph{Proc. CDC}, 2023, pp. 8320--8327.

\bibitem{Yang1993}
T.~H. Yang and E.~Polak, ``{Moving horizon control of nonlinear systems with
  input saturation, disturbances and plant uncertainty},'' \emph{Int. J.
  Control}, vol.~58, no.~4, pp. 875--903, 1993.

\bibitem{Alamir2017}
M.~Alamir, ``{Contraction-based nonlinear model predictive control formulation
  without stability-related terminal constraints},'' \emph{Automatica},
  vol.~75, pp. 288--292, 2017.

\bibitem{Worthmann2017}
K.~Worthmann, M.~W. Mehrez, G.~K. Mann, R.~G. Gosine, and J.~Pannek,
  ``{Interaction of open and closed loop control in MPC},'' \emph{Automatica},
  vol.~82, pp. 243--250, 2017.

\bibitem{Furnsinn2024}
A.~F\"urnsinn, C.~Ebenbauer, and B.~Gharesifard, ``{Flexible-step Model
  Predictive Control based on Generalized Lyapunov Functions},''
  \emph{Automatica}, vol. 175, p. 112215, 2025.

\bibitem{Furnsinn2023}
------, ``{Relaxed Feasibility and Stability Criteria for Flexible-Step MPC},''
  \emph{IEEE Control Syst. Lett.}, vol.~7, pp. 2851--2856, 2023.

\bibitem{Furnsinn2025}
------, ``{Flexible-step MPC for Switched Linear Systems with No Quadratic
  Common Lyapunov Function},'' \emph{IEEE Trans. Automat. Control}, 2025,
  submitted, see also https://arxiv.org/abs/2404.07870.

\bibitem{Feldbaum1960}
A.~A. Fel'dbaum, ``{Theory of dual control. Parts I, II, III, IV.}''
  \textit{Automation and Remote Control}, 21(9):874--880, 1960;
  21(11):1033--1039, 1960; 22(1):1--12, 1961; 22(2):109--121, 1961.

\bibitem{MichelBook}
A.~N. Michel, L.~Hou, and D.~Liu, \emph{{Stability of Dynamical Systems}},
  2nd~ed.\hskip 1em plus 0.5em minus 0.4em\relax Cham: Birkh\"auser, 2015.

\bibitem{Willems2005}
J.~C. Willems, P.~Rapisarda, I.~Markovsky, and B.~L. {De Moor}, ``{A note on
  persistency of excitation},'' \emph{Syst. Control Lett.}, vol.~54, no.~4, pp.
  325--329, 2005.

\bibitem{Yu2021}
Y.~Yu, S.~Talebi, H.~J. {van Waarde}, U.~Topcu, M.~Mesbahi, and
  B.~A\c{c}{\i}kme\c{s}e, ``{On Controllability and Persistency of Excitation
  in Data-Driven Control},'' in \emph{Proc. CDC}, 2021, pp. 6485--6490.

\bibitem{Zeng2023}
X.~Zeng, Z.~Liu, Z.~Du, N.~Ozay, and M.~Sznaier, ``On the hardness of learning
  to stabilize linear systems,'' \emph{2023 62nd IEEE Conference on Decision
  and Control (CDC)}, pp. 6622--6628, 2023.

\bibitem{Mityagin2020}
B.~S. Mityagin, ``{The Zero Set of a Real Analytic Function},''
  \emph{Mathematical Notes}, vol. 107, no.~3, pp. 529--530, 2020.

\end{thebibliography}

\appendix\hypertarget{sec:appendix}{}
In the subsequent lemma, we show that a large class of sequences of random vectors is almost surely persistently exciting (cf. \Cref{def:PE}). We include a proof of this statement since, surprisingly, we have not seen this elsewhere in the literature.
\begin{lemma}\label{lem:PE}
Let $\mathbb{U}$ be a finite-dimensional vector space. Let $T$ and $d$ be positive integers. Assume that there exists a persistently exciting sequence in $\mathbb{U}^{[\vphantom{]}0\colon\!\!T\vphantom{(})}$ of order $d$. Let $u_0,\ldots,u_{T-1}$ be independent $\mathbb{U}$-valued random variables on a probability space $(\Omega,\mathcal{F},\mathbb{P})$. Assume that, for every $i\in[\vphantom{]}0\colon\!\!T\vphantom{(})$, the probability distribution of $u_i$ is absolutely continuous with respect to the Lebesgue measure $\lambda$ on $\mathbb{U}$; that is, there exists a nonnegative Borel measurable function $\rho_i$ on $\mathbb{U}$ such that $\mathbb{P}[u_i\in\mathcal{B}]=\int_{\mathcal{B}}\rho_i\,\mathrm{d}\lambda$ for every Borel set $\mathcal{B}$ of $\mathbb{U}$. Then $u=(u_0,\ldots,u_{T-1})$ is almost surely persistently exciting of order $d$.
\end{lemma}
\begin{proof}
By the hypothesis, there exists a persistently exciting sequence $u^\ast\in\mathbb{U}^{[\vphantom{]}0\colon\!\!T\vphantom{(})}$ of order $d$. Notice that $md<{T-d}$, where $m$ denotes the dimension of $\mathbb{U}$. By the definition of persistency of excitation, there exist pairwise distinct elements $l(1)$, $\ldots$, $l(md)$ of $[0\colon\!\!T-d]$ such that the subsequences $u^\ast_{[\vphantom{]}l(1)\colon\!\!l(1)+d\vphantom{(})},\ldots,u^\ast_{[\vphantom{]}l(md)\colon\!\!l(md)+d\vphantom{(})}$ of $u^\ast$ form a basis of $\mathbb{U}^d$. As an abbreviation, we introduce the indexing set $K:=[\vphantom{]}l(1)\colon\!\!l(1)+d\vphantom{(})\cup\cdots\cup[\vphantom{]}l(md)\colon\!\!l(md)+d\vphantom{(})$ for those elements of $u^\ast$ that contribute to the above basis of $\mathbb{U}^d$. Define $v^\ast\in\mathbb{U}^K$ by $v^\ast_k:=u^\ast_k$ for every $k\in{K}$. For every $v\in\mathbb{U}^K$, let $L(v)$ be the endomorphism of $\mathbb{U}^d$ that assigns to the basis vector $u^\ast_{[\vphantom{]}l(i)\colon\!\!l(i)+d\vphantom{(})}$ the value $v_{[\vphantom{]}l(i)\colon\!\!l(i)+d\vphantom{(})}$ for every $i\in\{1,\ldots,md\}$. This defines a linear map $L$ from $\mathbb{U}^K$ to $\mathcal{L}(\mathbb{U}^d,\mathbb{U}^d)$. Define the real-valued function $p:=\det\circ{L}$ on $\mathbb{U}^K$. Since $p(v^\ast)=1$, the real-analytic function $p$ is not identically zero and therefore its zero set $\mathcal{N}\subseteq\mathbb{U}^K$ has Lebesgue measure zero (see, e.g., \cite{Mityagin2020}). Define a $\mathbb{U}^K$-valued random variable $v$ by $v_k:=u_k$ for every $k\in{K}$. Since the $\mathbb{U}$-valued random variables $v_k$ are assumed to be independent and since each of the $v_k$ is assumed to admit a density function $\rho_k$, also the $\mathbb{U}^K$-valued random variable $v$ admits a density function, which is denoted by $\rho$. Using that $\mathcal{N}$ has measure zero with respect to the Lebesgue measure $\mu$ on $\mathbb{U}^K$, it follows that
\begin{equation*}
\mathbb{P}[p(v)=0] \ = \ \mathbb{P}[v\in\mathcal{N}] \ = \ \int_\mathcal{N}\rho\,\mathrm{d}\mu \ = \ 0.
\end{equation*}
By the definition of $p$, this implies that $L(v)$ is almost surely invertible. Using the definitions of $L$ and $v$, we conclude that the subsequences $u_{[\vphantom{]}l(1)\colon\!\!l(1)+d\vphantom{(})},\ldots,u_{[\vphantom{]}l(md)\colon\!\!l(md)+d\vphantom{(})}$ of $u$ form almost surely a basis of $\mathbb{U}^d$. The claim follows.
\end{proof}
\end{document}